\documentclass[12pt]{article}
\usepackage{graphicx}
\usepackage{color}
\catcode`\@=11 \@addtoreset{equation}{section}

\catcode`\@=12

\newtheorem{Theorem}{Theorem}[section]
\newtheorem{Proposition}{Proposition}[section]
\newtheorem{Lemma}{Lemma}[section]
\newtheorem{Corollary}{Corollary}[section]
\newtheorem{Remark}{Remark}[section]
\newtheorem{Definition}{Definition}[section]

\newcommand{\bTheorem}[1]{
\begin{Theorem} \label{T#1} }
\newcommand{\eT}{\end{Theorem}}

\newcommand{\bProposition}[1]{
\begin{Proposition} \label{P#1}}
\newcommand{\eP}{\end{Proposition}}

\newcommand{\bLemma}[1]{
\begin{Lemma} \label{L#1} }
\newcommand{\eL}{\end{Lemma}}

\newcommand{\bCorollary}[1]{
\begin{Corollary} \label{C#1} }
\newcommand{\eC}{\end{Corollary}}

\newcommand{\bDefinition}[1]{
\begin{Definition} \label{D#1} }
\newcommand{\eD}{\end{Definition}}

\newcommand{\bRemark}[1]{
\begin{Remark} \label{R#1} }
\newcommand{\eR}{\end{Remark}}

\newcommand{\bFormula}[1]{
\begin{equation} \label{#1}}
\newcommand{\eF}{\end{equation}}

\newcommand{\Ov}[1]{\overline{#1}}

\newcommand{\DC}{C^\infty_c}

\newcommand{\vr}{\varrho}

\newcommand{\vt}{\vartheta}
\newcommand{\vu}{\vc{u}}
\newcommand{\vc}[1]{{\bf #1}}
\newcommand{\qed}{\bigskip \rightline {Q.E.D.} \bigskip}
\newcommand{\Div}{{\rm div}_x}
\newcommand{\Grad}{\nabla_x}

\newcommand{\tn}[1]{\mbox {\F #1}}
\newcommand{\dx}{{\rm d} {x}}
\newcommand{\dt}{{\rm d} t }

\newcommand{\dxdt}{\dx \ \dt}
\newcommand{\intO}[1]{\int_{\Omega} #1 \ \dx}

\newcommand{\bProof}{{\bf Proof: }}

\newcommand{\D}{{\cal D}}
\newcommand{\ep}{\varepsilon}

\font\F=msbm10 scaled 1000

\newcommand{\Del}{\Delta_x}

\definecolor{Cgrey}{rgb}{0.85,0.85,0.85}
\definecolor{Cblue}{rgb}{0.50,0.85,0.85}
\definecolor{Cred}{rgb}{1,0,0}
\definecolor{fancy}{rgb}{0.10,0.85,0.10}

\newcommand\Cbox[2]{%
    \newbox\contentbox%
    \newbox\bkgdbox%
    \setbox\contentbox\hbox to \hsize{%
        \vtop{
            \kern\columnsep
            \hbox to \hsize{%
                \kern\columnsep%
                \advance\hsize by -2\columnsep%
                \setlength{\textwidth}{\hsize}%
                \vbox{
                    \parskip=\baselineskip
                    \parindent=0bp
                    #2
                }%
                \kern\columnsep%
            }%
            \kern\columnsep%
        }%
    }%
    \setbox\bkgdbox\vbox{
        \color{#1}
        \hrule width  \wd\contentbox %
               height \ht\contentbox %
               depth  \dp\contentbox
        \color{black}
    }%
    \wd\bkgdbox=0bp%
    \vbox{\hbox to \hsize{\box\bkgdbox\box\contentbox}}%
    \vskip\baselineskip%
}

\date{}

\begin{document}


\title{Weak solutions to problems involving inviscid fluids}

\author{Eduard Feireisl \thanks{The research of E.F. leading to these results has received funding from the European Research Council under the European Union's Seventh Framework
Programme (FP7/2007-2013)/ ERC Grant Agreement 320078} }

\maketitle

\bigskip

\centerline{Institute of Mathematics of the Academy of Sciences of the Czech Republic}

\centerline{\v Zitn\' a 25, 115 67 Praha 1, Czech Republic}






\bigskip





\begin{abstract}

We consider an abstract functional-differential equation derived from the pressureless Euler system with variable coefficients that includes several
systems of partial differential equations arising in the fluid mechanics. Using the method of convex integration we show the existence of
infinitely many weak solutions for prescribed initial data and kinetic energy.

\end{abstract}



\section{Introduction}
\label{i}

The concept of \emph{weak solution} is indispensable in the mathematical theory of inviscid fluids, where solutions of the underlying non-linear systems of
partial differential equations are known to develop singularities in a finite lap of time no matter how smooth the initial data might be. The weak solutions
are being used even in the analysis of certain viscous fluids like the standard Navier-Stokes system, where a rigorous theory in the classical framework represents one of the major open problems of modern mathematics. In the absence of a sufficiently strong dissipative mechanism, solutions of non-linear systems of conservation laws may develop fast oscillations and/or concentrations that inevitably give rise to singularities of various types. As shown in the nowadays classical work of
Tartar \cite{T}, oscillations are involved in many problems, in particular in those arising in the context of inviscid fluids.

The well know deficiency of weak solutions is that they may not be uniquely determined in terms of the data and suitable admissibility criteria must be imposed in order to pick up
the physically relevant ones, cf. Dafermos \cite{D4}. Although most of the admissibility constraints are derived from fundamental physical principles
as the Second law of thermodynamics, their efficiency in eliminating the nonphysical solutions is still dubious, cf. Dafermos \cite{Daf4}.
Recently, DeLellis and Sz\' ekelyhidi \cite{DelSze3} developed the method previously known as \emph{convex integration} in the context
of fluid mechanics, in particular for the Euler system. Among other interesting results, they show the existence of infinitely many solutions to the
incompressible Euler system violating many of the standard admissibility criteria. Later, the method was adapted to the
compressible case by Chiodaroli \cite{Chiod}.

In this note, we introduce an abstract functional-differential equation that may be viewed as the pressureless Euler system with variable (functionally solution dependent) coefficients. We present an abstract version of the so-called oscillatory lemma and use it in order to show the existence of infinitely
many solutions adapting the method of \cite{DelSze3}. Various specific systems arising in fluid dynamics will be then identified as special cases of the abstract problem.

The paper is organized as follows. In Section \ref{E}, we introduce the abstract problem and formulate our main result proved in the remaining part of the paper. To this end, we adapt the apparatus of convex integration including the concept of \emph{subsolution} in Section \ref{C}. In Section \ref{O}, we present the
oscillatory lemma and show the existence of infinitely many solutions. Several specific examples  are discussed
in Section \ref{X}. Finally, Section \ref{A} addresses the problem of strong continuity of the weak solutions at the initial time.

\section{Abstract problem, main result}
\label{E}

The symbol $R^{N \times N}_{\rm sym}$ will denote the space of $N \times N$ symmetric matrices over the Euclidean space $R^N$, $N=2,3$,
$R^{N \times N}_{{\rm sym},0}$ is its subspace of those with zero trace. For two vectors $\vc{v}, \vc{w} \in R^N$, we denote
\[
\vc{v} \otimes \vc{w} \in R^{N \times N}_{\rm sym},\ [\vc{v} \otimes \vc{w}]_{i,j} = v_i v_j, \ \mbox{and}\
\vc{v}\odot \vc{w} \in R^{N \times N}_{{\rm sym},0}, \ \vc{v} \odot \vc{w} = \vc{v} \otimes \vc{w} - \frac{1}{N} \vc{v} \cdot \vc{w} \tn{I}.
\]

For the sake of simplicity, we suppose the physical space to be the ``flat'' torus
\[
\Omega = \left( [-1,1] \Big|_{\{-1, 1\}} \right)^N,
\]
meaning, the functions of $x \in \Omega$ are (2-)periodic in $R^N$.

\subsection{Abstract problem}

We consider the following problem:

Find a vector field $\vu \in C_{\rm weak}([0,T]; L^2(\Omega; R^N))$ satisfying
\bFormula{E1}
\partial_t \vc{u} + \Div \left( \frac{ (\vc{u} + \vc{h}[\vc{u}] ) \odot (\vc{u} + \vc{h}[\vc{u}] ) }{r[\vu]} + \tn{H}[\vu] \right) = 0,\
\Div \vu = 0
\ \mbox{in} \ \D'((0,T) \times \Omega; R^N),
\eF
\bFormula{E2}
\frac{1}{2} \frac{ | \vu + \vc{h} [\vu] |^2 }{r[\vu]}  (t, x) = e[\vu] (t,x) \ \mbox{for a.a.}\ (t,x) \in (0,T) \times \Omega,
\eF
\bFormula{E3}
\vu(0, \cdot) = \vu_0, \ \vu(T, \cdot) = \vu_T,
\eF
where $\vc{h}[\vu]$, $r[\vu]$, $\tn{H}[\vu]$, and $e[\vu]$ are given (nonlinear) operators.

\bRemark{AR1}

The problem (\ref{E1} -- \ref{E3}) is seemingly overdetermined as both the initial and the end state are prescribed. Moreover, the
associated ``kinetic energy'' is constrained by (\ref{E2}). Specific applications will be given in Section \ref{X}.

\eR

\bRemark{RE1}

The choice
\[
\vc{h} = 0, r = 1, \tn{H} = 0, e = e(t)
\]
gives rise to the pressureless (incompressible) Euler system
\[
\partial_t \vu + \Div (\vu \otimes \vu) = 0 , \ \Div \vu = 0
\]
with the prescribed kinetic energy
\[
\frac{1}{2} |\vu|^2 = e(t)
\]
studied by Chiodaroli \cite{Chiod} and  DeLellis, Sz\' ekelyhidi \cite{DelSze3}.

\eR

\bRemark{E2}

Note that a ``more complex'' problem
\bFormula{E5}
\partial_t \vc{u} + \Div \left( \frac{ (\vc{u} + \vc{h}[\vc{u}] ) \otimes (\vc{u} + \vc{h}[\vc{u}] ) }{r[\vu]} + \tn{H}[\vu] \right) +
\Grad \Pi [\vu] = 0,\
\Div \vu = 0
\eF
can be converted to (\ref{E1}), (\ref{E2}), with
\[
e[\vu] = Z[\vu](t) - \frac{N}{2} \Pi[\vu],
\]
where $Z$ is an arbitrary spatially homogeneous function.

\eR

\bRemark{E3}

The ``pressure'' $\Pi$ in (\ref{E5}) can be incorporated in $\tn{H}$ by solving the problem
\[
\Div \tn{H}_\Pi = \Grad \Pi \ \mbox{in}\ \Omega, \ \tn{H}_\Pi(x) \in R^{N\times N}_{{\rm sym},0},\ x \in \Omega.
\]
We can take, for instance, the solution of the Lam\' e system
\[
\tn{H}_\Pi = \Grad \vc{U} + \Grad \vc{U}^T - \frac{2}{N} \Div \vc{U} \tn{I}.
\]
As observed by Desvillettes and Villani \cite[Section IV.I, Proposition 11]{DesVil}, the vector field $\vc{U}$ is uniquely determined up to an additive constant. Of course, in order to preserve certain continuity of $\tn{H}_\Pi$, more regularity of $\Pi$ is needed.

\eR

The quantities $\vc{h}$, $r$, $\tn{H}$, and $e$ are operators depending on the solution $\vu$. In order to specify their
properties, we introduce the following definition:

\bDefinition{R1}

Let $Q \subset (0,T) \times \Omega$ be an open set such that
\[
|Q| = |(0,T) \times \Omega|.
\]
An operator
\[
b: C_{\rm weak}([0,T]; L^2(\Omega; R^N)) \cap L^\infty((0,T) \times \Omega; R^N) \to C_b (Q, R^M)
\]
is \emph{$Q-$continuous} if:

\begin{itemize}

\item $b$ maps bounded sets in $L^\infty((0,T) \times \Omega; R^N)$ on bounded sets in $C_b (Q, R^M)$;

\item $b$ is continuous, specifically,

\bFormula{E4a}
\begin{array}{c}
b[\vc{v}_n] \to b[\vc{v}] \ \mbox{in} \ C_b(Q; R^M) \ \mbox{(uniformly for $(t,x) \in Q$ ) } \\ \\
\mbox{whenever}\\ \\
\vc{v}_n \to \vc{v} \ \mbox{in}\ C_{\rm weak}([0,T]; L^2(\Omega; R^N)) \ \mbox{and weakly-(*) in} \ L^\infty((0,T) \times \Omega; R^N);
\end{array}
\eF

\item $b$ is causal (non-anticipative), meaning
\bFormula{E4b}
\vc{v}(t, \cdot) = \vc{w}(t,\cdot) \ \mbox{for}\ 0 \leq t \leq \tau \leq T
\ \mbox{implies} \ b[\vc{v}] = b[\vc{w}] \ \mbox{in}\  \left[ (0, \tau] \times \Omega \right] \cap Q.
\eF

\end{itemize}

\eD

In this paper, we suppose
\bFormula{E4}
\begin{array}{c}
\vc{h} = \vc{h}[\vu]: C_{\rm weak}([0,T]; L^2(\Omega; R^N)) \to C_b
(Q; R^N ),\\ \\
r = r[\vu] : C_{\rm weak}([0,T]; L^2(\Omega; R^N)) \to C_b (Q; R),\ r > 0,\\ \\
e = e[\vu] : C_{\rm weak}([0,T]; L^2(\Omega; R^N)) \to C_b(Q; R),\ e \geq 0, \\ \\
\tn{H} = \tn{H}[\vu]: C_{\rm weak}([0,T]; L^2(\Omega; R^N)) \to  C_b(Q; R^{N \times N}_{{\rm sym},0} )
\end{array}
\eF
are given $Q-$continuous operators for a certain open set $Q$.

\subsection{Subsolutions}

Before stating our main result concerning solvability of problem (\ref{E1}--\ref{E3}), it is convenient to introduce the
set of \emph{subsolutions}.
Let $\lambda_{\rm max}[\tn{A}]$ denote the maximal eigenvalue of a matrix $\tn{A} \in R^{N \times N}_{\rm sym}$.
Similarly to DeLellis and Sz\' ekelyhidi \cite{DelSze3}, we introduce the set of \emph{subsolutions}:
\bFormula{E6}
X_0 = \left\{ \vc{v} \ \Big| \ {\vc{v}} \in C_{\rm weak}([0,T]; L^2(\Omega; R^N)) \cap L^\infty((0,T) \times \Omega;R^N) \vphantom{\frac{1}{2}},\
\vc{v}(0, \cdot) = \vu_0 , \ \vc{v}(T, \cdot) = \vu_T, \right.
\eF
\[
\partial_t \vc{v} + \Div \tn{F} = 0,
 \ \Div \vc{v} = 0 \ \mbox{in}\ \D'((0,T) \times \Omega; R^N),
\ \mbox{for some}\ \tn{F} \in L^\infty((0,T) \times \Omega; R^{N \times N}_{{\rm sym},0} ),
\]
\[
\vc{v} \in C(Q; R^N),\ \tn{F}  \in C(Q; R^{N \times N}_{{\rm sym},0} ),
\]
\[
\sup_{(t,x) \in Q, t > \tau}
\left.
\frac{N}{2}\lambda_{\rm max} \left[ \frac{ (\vc{v} + \vc{h}[\vc{v}] ) \otimes (\vc{v} + \vc{h}[\vc{v}] ) }{r[\vc{v}]} - \tn{F} + \tn{H}[\vc{v}] \right]
- e[\vc{v}] < 0 \ \mbox{for any}\ 0 < \tau < T \vphantom{\frac{1}{2}}  \right\}.
\]

\bRemark{E5}

Note that, in contrast with \cite{DelSze3},  the inequality
\[
\frac{N}{2}\lambda_{\rm max} \left[ \frac{ (\vc{v} + h[\vc{v}] ) \otimes (\vc{v} + h[\vc{v}] ) }{r[\vc{v}]} - \tn{F} + \tn{H}[\vc{v}] \right]
 < e[\vc{v}]
\]
is satisfied only on the open set $Q$, where all quantities are continuous. Moreover, the inequality is strict on any open
time interval $(\tau, T)$, $0 < \tau < T$.

\eR

\subsection{Main result}
\label{MR}

We are ready to state our main result.

\bTheorem{E1}

Let the operators $\vc{h}$, $r$, $\tn{H}$, and $e$ given by (\ref{E4}) be $Q-$continuous, where $Q \subset \left[ (0,T) \times \Omega \right]$ is an open set,
\[
|Q| = |(0,T) \times \Omega|.
\]
In addition, suppose that $r[\vc{v}] > 0$ and that the mapping $\vc{v} \mapsto 1/r[\vc{v}]$ is continuous in the sense specified in (\ref{E4a}).
Finally, assume that the set of subsolutions $X_0$ is \emph{non-empty} and \emph{bounded} in $L^\infty((0,T) \times \Omega; R^N)$.

Then problem (\ref{E1} -- \ref{E3}) admits infinitely many solutions.

\eT

The next two sections will be devoted to the proof of Theorem \ref{TE1}. For the set of subsolutions to be non-empty, the energy $e$ must be chosen
large enough. For instance, taking $\vc{u}_0 = \vu_T \in C(\Omega; R^N)$, $\Div \vu_0 = 0$ we check easily that $X_0$ is non-empty, specifically
$\vu_0 \in X_0$, as soon as
\bFormula{E7}
\frac{N}{2}\lambda_{\rm max} \left[ \frac{ (\vu_0 + h[\vu_0] ) \otimes (\vu_0 + h[\vu_0] ) }{r[\vu_0]} + \tn{H}[\vu_0] \right]
 < e[\vu_0].
\eF

Recalling the purely algebraic inequality (cf. \cite{DelSze3})
\bFormula{E7a}
\frac{1}{2} \frac{ |\tilde{\vc{h}} |^2 }{\tilde r} \leq
\frac{N}{2}\lambda_{\rm max} \left[ \frac{ \tilde{\vc{h}}  \otimes \tilde{\vc{h}} }{\tilde r} - \tilde{H} \right],
\eF
where the equality holds only if
\bFormula{E7b}
\tilde{H} = \frac{ \tilde{\vc{h}}  \otimes \tilde{\vc{h}} }  {\tilde r} - \frac{1}{N} \frac{ |\tilde{\vc{h}} |^2 } {\tilde r} \tn{I},
\eF
we get from (\ref{E7}) that
\[
\frac{1}{2} \frac{ |\vu_0 + h [\vu_0] |^2 }{r[\vu_0]} < e[\vu_0],
\]
meaning the relation (\ref{E2}) is \emph{violated} at the initial time. This is the undesirable initial
``energy jump'' characteristic for the weak solutions obtained by the method of convex integration. A possible remedy for this problem will be discussed in  Section \ref{A}.

\section{Convex integration}
\label{C}

As the set $X_0$ is bounded, there exists a positive constant $\Ov{e}$ such that
\bFormula{C1-}
e[\vc{v}] \leq \Ov{e} \ \mbox{for any}\ \vc{v} \in X_0.
\eF
Under the hypotheses of Theorem \ref{TE1} we may define a topological space $\Ov{X}_0$
as the closure of the space of subsolutions $X_0$ with respect to the (metrizable) topology of $C_{\rm weak}([0,T]; L^2(\Omega; R^N))$. Accordingly,
$\Ov{X}_0$ is a (non-empty) complete metric space with the distance of two functions $\vc{v}$, $\vc{w}$ given by
\[
\sup_{t \in [0,T]} d[\vc{v}(t, \cdot); \vc{w}(t,\cdot)],
\]
where $d$ is the metrics induced by the weak topology on bounded sets of the Hilbert space $L^2(\Omega; R^N)$.
Note that, in view of (\ref{E7a}), (\ref{C1-}) and boundedness of all operators involved in the definition of $X_0$, the associated fluxes
$\tn{F}$ are bounded in $L^\infty$, in particular,
\bFormula{EQ1}
\partial_t \vc{v} + \Div \tn{F} = 0,\ \Div \vc{v} = 0 \ \mbox{in}\ \D'((0,T) \times \Omega; R^N),
\eF
for any $\vc{v} \in \Ov{X}_0$, where the flux $\tn{F} \in L^\infty((0,T) \times \Omega; R^{N\times N}_{{\rm sym},0})$ can be obtained as a weak limit of fluxes in $X_0$. Moreover, by convexity of the function
\[
\frac{N}{2}\lambda_{\rm max} \left[ \frac{ (\vc{v} + \vc{h} ) \otimes (\vc{v} + \vc{h} ) }{r} - \tn{F} + \tn{H} \right]
\]
in $\vc{v}$ and $\tn{F}$, we get
\bFormula{EQ2}
\frac{N}{2}\lambda_{\rm max} \left[ \frac{ (\vc{v} + \vc{h}[\vc{v}] ) \otimes (\vc{v} + \vc{h}[\vc{v}] ) }{r[\vc{v}]} - \tn{F} + \tn{H}[\vc{v}] \right] \leq e[\vc{v}] \ \mbox{a.a. in}\ (0,T) \times \Omega.
\eF

Next,
we introduce a countable family of functionals
\[
I_n [\vc{v}] = \int_{1/n}^{T} \intO{ \left[ \frac{1}{2} \frac{ |\vc{v} + \vc{h}[\vc{v}]  |^2 }{r[\vc{v}]} - e[\vc{v}] \right] } : \Ov{X}_0 \to (-\infty,0].
\]
In accordance with the hypotheses (\ref{E4}), each $I_n$ can be seen as a lower semi-continuous functional on $\Ov{X}_0$. In particular, by means of
Baire's category argument, the set
\[
\mathcal{S} = \cap_{n > 0} \left\{ \vc{v} \in \Ov{X}_0 \ \Big|\ \vc{v} \ \mbox{is a point of continuity of} \ I_n \right\}
\]
has infinite cardinality.

In the next section, we show that
\bFormula{C1}
\mbox{if}\
\vc{v} \ \mbox{is a point of continuity of}\ I_n \ \mbox{in}\ \Ov{X}_0, \ \mbox{then} \ I_n[\vc{v}] = 0.
\eF
In accordance with (\ref{E7a}), (\ref{E7b}), combined with the previous observations stated in
(\ref{EQ1}), (\ref{EQ2}), this implies that $\mathcal{S}$ consists of weak solutions to problem (\ref{E1}--\ref{E3}).
Consequently, the proof of Theorem \ref{TE1} reduces to showing (\ref{C1}).

\section{Oscillatory lemma, infinitely many solutions}
\label{O}

In accordance with the previous discussion,
the final step in the proof of Theorem \ref{TE1} is to show (\ref{C1}).
The main tool we shall use is the following variant of the oscillatory lemma (cf. De Lellis
and Sz\' ekelyhidi \cite[Proposition 3]{DelSze3}, Chiodaroli \cite[Section 6, formula (6.9)]{Chiod}) proved in \cite[Lemma 3.1]{DoFeMa} :

\bLemma{O1}

Let $U \subset R \times R^N$, $N=2,3$ be a bounded open set. Suppose that
\[
\tilde {\vc{h}} \in C(U; R^N), \ \tilde {\tn{H}} \in C(U; R^{N \times N}_{{\rm sym},0}), \ \tilde e, \ \tilde r \in C(U),\ \tilde r > 0, \ \tilde e \leq \Ov{e} \ \mbox{in}\ U
\]
are given such that
\bFormula{O12-}
\frac{N}{2} \lambda_{\rm max} \left[ \frac{ \tilde {\vc{h}} \otimes \tilde{\vc{h}} }{\tilde{r}} - \tilde{\tn{H}} \right] < \tilde e  \ \mbox{in}\
U.
\eF

Then there exist sequences
\[
\vc{w}_n \in \DC (U;R^N), \ \tn{G}_n \in \DC(U; R^{N \times N}_{\rm sym,0}), \ n = 0,1,\dots
\]
such that
\[
\partial_t \vc{w}_n + \Div \tn{G}_n = 0 , \ \Div \vc{w}_n = 0 \ \mbox{in}\ R^N,
\]
\bFormula{O1-}
\frac{N}{2} \lambda_{\rm max} \left[ \frac{ (\tilde {\vc{h} } + \vc{w}_n) \otimes ( \tilde{\vc{h}} + \vc{w}_n)}{\tilde{r}} - (\tilde {\tn{H}} + \tn{G}_n) \right] < \tilde{e}  \ \mbox{in}\ U,
\eF
and
\bFormula{O1}
\vc{w}_n \to 0 \ \mbox{weakly in}\ L^2(U; R^N),\
\liminf_{n \to \infty} \int_{U} \frac{ | \vc{w}_n |^2 }{ \tilde{r}} \ \dxdt \geq \Lambda(\Ov{e}) \int_{U} \left( \tilde{e} - \frac{1}{2} \frac{| \tilde{\vc{h}} |^2}{
\tilde r} \right)^2 \ \dxdt
\eF
for a certain $\Lambda (\Ov{e}) > 0$ depending only on the energy upper bound $\Ov{e}$.

\eL

\bRemark{AR5}

Note that Lemma \ref{LO1} applies to continuous, not necessarily bounded, functions on the open set $U$.

\eR

With Lemma \ref{LO1} at hand, we may show the following result that contains (\ref{C1}) as a particular case.

\bLemma{O2}
Let
\[
I_D= \int_D \left[ \frac{1}{2} \frac{ |\vc{v} + \vc{h}[\vc{v}]  |^2 }{r[\vc{v}]} - e[\vc{v}] \right] \dxdt : \Ov{X}_0 \to (-\infty,0]
\]
be a functional defined on an open set $D \subset \left[ (\tau, T) \times \Omega \right] \cap Q$, $0 < \tau < T$.

Then $I_D$ vanishes at any of its points of continuity.

\eL

\bProof

Arguing by contradiction we assume that $\vc{v} \in \Ov{X}_0$ is a point of continuity of $I_D$ such that
\[
I_D[\vc{v}] < 0.
\]
Since $I_D$ is continuous at $\vc{v}$, there is a sequence $\{ \vc{v}_m \}_{m=1}^\infty \subset X_0$ (with the associated fluxes $\tn{F}_m$) such that
\[
\vc{v}_m \to \vc{v} \ \mbox{in}\ C_{\rm weak}([0,T]; L^2(\Omega;R^N)), \ I_D[\vc{v}_m] \to I_D[\vc{v}] \ \mbox{as}\ m \to \infty.
\]

As $[\vc{v}_m, \tn{F}_m]$ are subsolutions and $\tau > 0$,  we get, thanks to (\ref{E6}),
\[
\frac{N}{2}\lambda_{\rm max} \left[ \frac{ (\vc{v}_m + \vc{h}[\vc{v}_m]) \otimes (\vc{v}_m + \vc{h}[\vc{v}_m]) }{r[\vc{v}_m]} - \tn{F}_m
+ \tn{H}[\vc{v}_m] \right]
\]
\[
< e[\vc{v}_m] - \delta_m \ \mbox{in}\ \left[ [\tau, T) \times \Omega) \right] \cap Q\  \mbox{for some}\ \delta_m \searrow  0.
\]

Now,
fixing $m$ for a while, we apply Lemma \ref{LO1} with
\[
N=2,3, \ U = D,\ \tilde r = r[\vc{v}_m], \ \tilde {\vc{h}} = \vc{v}_m + \vc{h}[\vc{v}_m], \ \tilde{\tn{H}} = \tn{F}_m - \tn{H}[\vc{v}_m],\
 \tilde{e} = e[\vc{v}_m] - \delta_m.
\]
For $\left\{ [\vc{w}_{m,n}, \tn{G}_{m,n}] \right\}_{n=1}^\infty$ the quantities resulting from the conclusion of Lemma \ref{LO1}, we set
\[
\vc{v}_{m,n} = \vc{v}_m + \vc{w}_{m,n}, \ \tn{F}_{m,n} = \tn{F}_m + \tn{G}_{m,n}.
\]

Obviously,
\[
\partial_t \vc{v}_{m,n} + \Div \tn{F}_{m,n} = 0,\ \Div \vc{v}_{m,n} = 0 \ \mbox{in}\ \D'((0,T) \times \Omega),\ \vc{v}_{m,n} (0, \cdot) = \vc{v}_0,
\ \vc{v}_{m,n}(T, \cdot) = \vc{v}_T.
\]

Moreover, in accordance with (\ref{O1-}) and the fact that $\vc{w}_n$, $\tn{G}_{m,n}$ vanish outside $D$,
\[
\frac{N}{2}\lambda_{\rm max} \left[ \frac{ (\vc{v}_{m,n} + \vc{h}[\vc{v}_m]) \otimes (\vc{v}_{m,n} + \vc{h}[\vc{v}_m] ) }{r[\vc{v}_m]} - \tn{F}_{m,n}
+ \tn{H}[\vc{v}_m] \right]
< {e}[\vc{v}_m] - \delta_m\ \mbox{in}\ [\tau, T) \times \Omega \cap Q,
\]
and, by virtue of the causality property (\ref{E4b}),
\[
\sup_{(t,x) \in Q, s < t \leq \tau}
\frac{N}{2}\lambda_{\rm max} \left[ \frac{ (\vc{v}_{m,n} + \vc{h}[\vc{v}_{m,n}]) \otimes (\vc{v}_{m,n} + \vc{h}[\vc{v}_{m,n}] ) }{r[\vc{v}_{m,n}]} - \tn{F}_{m,n}
+ \tn{H}[\vc{v}_{m,n}] \right] - {e}[\vc{v}_{m,n}] < 0
\]
for any $0 < s< \tau$.
Consequently, in view of continuity of the operators $\vc{v} \mapsto \vc{h}[\vc{v}], r[\vc{v}], e[\vc{v}], \tn{H}[\vc{v}]$ specified in (\ref{E4a}),
we may infer that
for each $m$ there exists $n=n(m)$ such that
\[
\vc{v}_{m,n(m)} \in X_0, \ m=1,2,\dots
\]
Moreover, by virtue of (\ref{O1}), we may suppose that
\[
\vc{v}_{m,n(m)} \to \vc{v} \ \mbox{in}\ C_{\rm weak}([0,T]; L^2(\Omega;R^2))
\]
in particular,
\bFormula{O2}
I_D[\vc{v}_{m,n(m)} ] \to I_D[\vc{v}]
\eF
as $m \to \infty$.

Finally, using again the conclusion of Lemma \ref{LO1} combined with Jensen's inequality, we observe that the sequence $\vc{v}_{m,n(m)}$ can be taken in such a way that
\[
\liminf_{m \to \infty} I_D [\vc{v}_{m,n(m)}] = \liminf_{m \to \infty} \int_D \left( \frac{1}{2} \frac{ |\vc{v}_m + \vc{w}_{m,n(m)} +
\vc{h}[\vc{v}_m + \vc{w}_{m,n(m)}]|^2 }{r[\vc{v}_m + \vc{w}_{m,n(m)}]} - e [\vc{v}_m + \vc{w}_{m,n(m)}] \right) \dxdt
\]
\[
= \lim_{m \to \infty} \int_D \left( \frac{1}{2} \frac{ |\vc{v}_m + \vc{h}[\vc{v}_m + \vc{w}_{m,n(m)} ]|^2 }{r [\vc{v}_m + \vc{w}_{m,n(m)}]} - e [\vc{v}_m + \vc{w}_{m,n(m)}] \right) \dxdt
\]
\[
 + \liminf_{m \to \infty} \int_D \frac{1}{2} \frac{ |\vc{w}_{m,n(m)}|^2 }{r [\vc{v}_m + \vc{w}_{m,n(m)}]} \dxdt
\]
\[
\geq I_D [\vc{v}] + \frac{\Lambda (\Ov{e})}{2} \liminf_{m \to \infty} \int_D \left( e([\vc{v}_m])  -\delta_m - \frac{1}{2} \frac{|\vc{v}_m + \vc{h}[\vc{v}_m]|^2 }
{r[\vc{v}_m] } \right)^2  \dxdt
\]
\[
\geq
I_D[\vc{v}] + \frac{\Lambda (\Ov{e})}{2 |G|} \liminf_{m \to \infty} \left( \int_D \left( e([\vc{v}_m])  -\delta_m - \frac{1}{2} \frac{|\vc{v}_m + \vc{h}[\vc{v}_m] |^2 }{r[\vc{v}_m] } \right) \dxdt   \right)^2
=
I_D[\vc{v}] +  \frac{\Lambda (\Ov{e})}{2 |G|} \left( I_D[\vc{v}] \right)^2,
\]
which is compatible with (\ref{O2}) only if $I_D[\vc{v}] = 0$.

\rightline{\qed}

We have shown (\ref{C1}); whence Theorem \ref{TE1}.

\section{Examples}
\label{X}

There are many systems arising in mathematical fluid dynamics that can be written in the abstract form (\ref{E1} -- \ref{E3}). We review some of them already studied in the available literature.

\subsection{Euler-Fourier system}
\label{EF}

The Euler-Fourier system describes the time evolutions of the mass density $\vr = \vr(t,x)$, the velocity $\vu = \vu(t,x)$, and the (absolute) temperature
$\vt = \vt(t,x)$:
\bFormula{X1}
\partial_t \vr + \Div (\vr \vu) = 0,
\eF
\bFormula{X2}
\partial_t (\vr \vu) + \Div (\vr \vu \otimes \vu) + \Grad (\vr \vt) = 0,
\eF
\bFormula{X3}
\frac{3}{2} \Big( \partial_t (\vr \vt) + \Div (\vr \vt \vu) \Big) - \Delta \vt = - \vr \vt \Div \vu.
\eF

Following \cite{ChiFeiKre} we first write the momentum $\vr \vu$ as its Helmholtz decomposition
\[
\vr \vu = \vc{v} + \Grad \Phi, \ \Div \vc{v} = 0.
\]
Accordingly, we may \emph{fix} the density $\vr$ and the acoustic potential $\Phi$ so that
\[
\partial_t \vr + \Delta \Phi = 0 \ \mbox{holds,}
\]
meaning equation (\ref{X1}) is satisfies as $\Div \vc{v} = 0$.

With $\vr$, $\Phi$ given we may determine the temperature field $\vt = \vt[\vc{v}]$ as the (unique solution) of (\ref{X3}), specifically
\[
\frac{3}{2} \Big( \vr \partial_t \vt + (\vc{v} + \Grad \Phi) \cdot \Grad \vt \Big) - \Delta \vt = - \vr \vt \Div \left( \frac{1}{\vr} \left( \vc{v} + \Grad \Phi \right) \right)
\]
endowed with appropriate initial data.

Finally, we rewrite (\ref{X2}) in the form
\bFormula{X4}
\partial_t \vc{v} + \Div\left( \frac{ (\vc{v} + \Grad \Phi) \otimes (\vc{v} + \Grad \Phi ) }{\vr} \right) + \Grad \left( \partial_t \Grad \Phi +
\vr \vt[\vc{v}] \right) = 0.
\eF
Fixing the ``energy'' so that
\bFormula{X5}
\frac{1}{2} \frac{ |\vc{v} + \Grad \Phi |^2 }{r} = e[\vc{v}] \equiv Z - \frac{N}{2} \left( \partial_t \Grad \Phi - \vr \vt [\vc{v}] \right),
\eF
where $Z = Z(t)$ is a suitable spatially homogeneous function, we reduce (\ref{X4}) to
\bFormula{X6}
\partial_t \vc{v} + \Div\left( \frac{ (\vc{v} + \Grad \Phi) \odot (\vc{v} + \Grad \Phi ) }{\vr} \right) = 0, \ \Div \vc{v} = 0,
\eF
which is an equation in the form (\ref{E1}).

With certain effort, it is possible to show that the hypotheses of Theorem \ref{TE1} are satisfied for $Q= (0,T) \times \Omega$,
and we obtain the following result, see \cite[Theorem 3.1]{ChiFeiKre}:

\bTheorem{X1}
Let $T > 0$ be given, along with the initial data
\bFormula{X7}
\vr(0, \cdot) = \vr_0 \in C^3(\Omega), \ \vr_0 > 0, \ \vt(0,\cdot) = \vt_0 \in C^2(\Omega), \ \vt_0 > 0, \ \vu(0, \cdot) = \vu_0 \in C^3(\Omega; R^N),
\eF
\[
\Omega = \left( [-1; 1]_{\{ -1; 1 \} } \right)^N, \ N=2,3.
\]

Then the Euler-Fourier system (\ref{X1} -- \ref{X3}) admits infinitely many weak solutions in $(0,T) \times \Omega$ emanating from the same initial data
(\ref{X7}).
\eT

As already pointed out, the solutions obtained in Theorem \ref{TX1} may be non-physical in the sense they violate the principle of energy conservation. However,
this drawback can be removed at least for certain initial data. We will discuss this issue in Section \ref{KON}.

\subsection{Quantum fluids}

The Euler-Korteweg-Poisson system describes the time evolution of the density $\vr = \vr(t,x)$ and the momentum $\vc{J} = \vc{J}(t,x)$ of
an inviscid fluid:

\bFormula{X8}
\partial_t \vr + \Div \vc{J} = 0,
\eF
\bFormula{X9}
\partial_t \vc{J} + \Div \left( \frac{\vc{J} \times \vc{J} }{\vr} \right) + \Grad p(\vr) =  - \alpha \vc{J} +  \vr \Grad \left(
K(\vr) \Del \vr + \frac{1}{2} K'(\vr) |\Grad \vr |^2 \right) + \vr \Grad V,
\eF
\bFormula{X10}
\Del V = \vr - \Ov{\vr},
\eF
where $K : (0, \infty) \to (0, \infty)$ is a given function, see Audiard \cite{Audi2}, Benzoni-Gavage et al. \cite{Benz2}, \cite{Benz1}. The choice
$K = \Ov{K} > 0$ yields the standard equations of an inviscid capillary fluid (see Bresch et al.
\cite{BDD}, Kotchote \cite{Kot2}, \cite{Kot1}), while $K(\vr) = \frac{\hbar}{4 \vr}$ gives rise to the quantum fluid system (see for instance Antonelli and Marcati \cite{AntMar1}, \cite{AntMar2}, J\" ungel \cite[Chapter 14]{Jun} and the references therein).

For
\[
\chi(\vr) = \vr K(\vr),
\]
it can be shown that system (\ref{X8} -- \ref{X10}) can be recast in the form
\bFormula{X11}
\partial_t  \vc{v}  + \Div \left( \frac{ ( \vc{v} + \vc{h} ) \otimes (\vc{v} + \vc{h}) }{ r } + \tn{H} \right)
 + \Grad \Pi = 0,
\eF
with
\[
r = e^t \vr, \ \vc{h} = e^t \Grad M,
\]
\[
\tn{H}(t,x) = 4 e^t \left( \chi(\vr) \Grad \sqrt{\vr} \otimes \Grad \sqrt{\vr} - \frac{1}{3} \chi(\vr) |\Grad \sqrt{\vr} |^2 \tn{I} -
\frac{1}{4}\Grad V \otimes \Grad V + \frac{1}{12} |\Grad V|^2 \tn{I} \right),
\]
and
\[
\Pi(t,x) = e^t \left( p(\vr) + \partial_t M + M - \chi(\vr) \Del \vr - \frac{1}{2} \chi'(\vr) |\Grad \vr|^2 + \frac{4}{3}
\chi(\vr) |\Grad \sqrt{\vr} |^2 - \Ov{\vr} V +\frac{1}{6} |\Grad V|^2 \right),
\]
where $\vr$ and $M$ are suitably chosen functions, see \cite{DoFeMa}.

Now, Theorem \ref{TE1} can be applied to obtain the following result, see \cite[Theorem 2.1]{DoFeMa} and the proof therein.

\bTheorem{X2}
Let $T > 0$ be given. Suppose that $p$ and $\chi$ satisfy
\[
p \in C^1[0, \infty) \cap C^2(0, \infty), \ p(0) = 0, \ \chi \in C^2[0, \infty), \ \chi > 0 \ \mbox{in}\ (0, \infty).
\]
Let the initial data be given such that
\bFormula{X12}
\vr(0, \cdot) = \vr_0 = r_0^2, \ r_0 \in C^2(\Omega), \ {\rm meas} \left\{ x \in \Omega \ \Big| \ r_0(x) = 0 \right\} = 0,
\eF
\bFormula{X13}
\vc{J}(0, \cdot) =
\vc{J}_0 = \vr_0 \vc{U}_0, \ \vc{U}_0 \in C^3( \Omega; R^3).
\eF

Then the initial value problem (\ref{X8}--\ref{X10}), (\ref{X12}), (\ref{X13}) admits infinitely many weak solutions in $(0,T) \times \Omega$.

\eT

In the situation described in Theorem \ref{TX2}, the set $Q$ must be taken
\[
Q = (0,T) \times \Omega \setminus \left\{ (t,x) \ \Big| \ \vr(t,x) = 0 \right\}.
\]

\subsection{Binary mixtures of compressible fluids}

We consider a physically motivated regularization of the Euler equations proposed in the seminal paper  by Lowengrub and Truskinovsky
\cite{LoTr}. The model describes the motion of a mixture of two immiscible compressible fluids in terms of the density $\vr = \vr(t,x)$, the
velocity $\vu = \vu(t,x)$, and the concentration difference $c = c(t,x)$. The fluid is described by means of the standard Euler system coupled with the Cahn-Hilliard equation describing the evolution of $c$:

\bFormula{X14}
\partial_t \vr + \Div (\vr \vu) = 0,
\eF
\bFormula{X15}
\partial_t(\vr \vu) + \Div (\vr \vu \otimes \vu) + \Grad p_0(\vr,c) = \Div \left( \vr \Grad c \otimes \Grad c - \frac{\vr}{2} |\Grad c|^2 \tn{I} \right),
\eF
\bFormula{X16}
\partial_t (\vr c) + \Div (\vr c \vu) = \Delta \left( \mu_0(\vr, c) - \frac{1}{\vr} \Div \left( \vr \Grad c \right) \right),
\eF
where
\bFormula{X17}
p_0 (\vr,c) = \vr^2 \frac{\partial f_0(\vr, c) }{\partial \vr} ,\ \mu_0 (\vr, c) = \frac{\partial f_0 (\vr ,c) }{\partial c}
\eF
for a given free energy function $f_0$.
The system is neither purely hyperbolic nor parabolic as the dissipation mechanism acts in a very subtle way through the coupling of the Euler and the Cahn-Hilliard systems.

The machinery of convex integration can be applied, first fixing $\vr$ and $\Phi$, similarly to Section \ref{EF},
to solve
\[
\partial_t \vr + \Delta \Phi = 0,
\]
then taking $c = c[\vc{v}]$, $\Div \vc{v} = 0$ to be the unique solution of the equation
\[
\partial_t (\vr c) + \Div (\vr c (\vc{v} + \Grad \Phi) ) = \Delta \left( \mu_0(\vr, c) - \frac{1}{\vr} \Div \left( \vr \Grad c \right) \right).
\]

Accordingly, we obtain
\bFormula{X18}
\vc{v}(0, \cdot) = \vc{v}_0, \ \Div \vc{v} = 0 ,
\eF
\bFormula{X19}
\partial_t  \vc{v}  + \Div \left( \frac{ (\vc{v} + \Grad \Phi) \odot (\vc{v} + \Grad \Phi) }{\vr} - \vr \left( \Grad c[\vc{v}] \odot \Grad c[\vc{v}] \right)
\right) = 0,
\eF
\bFormula{X20}
\frac{1}{2} \frac{ |\vc{v} + \Grad \Phi|^2 }{\vr}   = e [\vc{v}] \equiv Z(t) - \frac{N}{2} \left( \frac{1}{6} |\Grad c[\vc{v}] |^2 + p_0(\vr, c[\vc{v}] ) + \partial_t \Grad \Phi
 \right) ,
\eF
where $Z$ is a spatially homogeneous function.

Theorem \ref{TE1} yields the following result, see \cite{EF200} for details:

\bTheorem{X3}
Let the potential $f_0 = f_0(\vr, c)$ satisfy
\[
f_0(\vr,c) = H(c) + \log(\vr) \left( \alpha_1 \frac{1 - c}{2} + \alpha_2 \frac{1 + c}{2} \right),\
H \in C^2(R),\ |H''(c) | \leq \overline{H} \ \mbox{for all} \ c \in R^1.
\]

Then for any choice of initial conditions
\[
\vr(0, \cdot) = \vr_0 \in C^3(\Omega), \ \inf_{\Omega} \vr_0 > 0,\ \vc{u}(0, \cdot) = \vu_0 \in C^3(\Omega;R^3), \
c(0, \cdot) = c_0 \in C^2(\Omega),
\]
the problem (\ref{X14} - \ref{X16}) admits infinitely many weak solutions in $(0,T) \times \Omega$.

\eT

\section{Continuity at the initial time, admissible solutions}
\label{A}

The major drawback of the construction delineated in the previous part of the paper and the main reason why the weak solution obtained via
convex integration can be eliminated as physically unacceptable is the energy jump at the initial time discussed in Section \ref{MR}.
On the other hand, however, once a subsolution $\vc{v}$ along with the associated energy $e[\vc{v}]$ are obtained, it is possible to show
the existence of another subsolution defined on a possibly shorter time interval for which the initial energy is attained. Such a subsolution can be then used in the process of convex integration to produce weak solutions that are \emph{strongly continuous} at the initial time
and dissipate energy.

We first state the result for the abstract system and then shortly comment on possible applications.
Modifying slightly the procedure used in the proof of Theorem \ref{TE1} we can show the following assertion:

\bTheorem{A1}

In addition to the hypotheses of Theorem \ref{TE1}, suppose that
\bFormula{HYP}
\left| \left\{ x \in \Omega \ \Big| (t,x) \in Q \right\} \right| = |\Omega| \ \mbox{for any}\ 0 < t < T.
\eF

Then there exists a set of times $\mathcal{R} \subset (0,T)$ dense in $(0,T)$ such that for
any $\tau \in \mathcal{R}$ there is $\vc{v} \in \Ov{X}_0$ with the following properties:
\begin{itemize}
\item
\bFormula{pr1}
\vc{v} \in C_b(\left[ (0,\tau) \cup (\tau, T) \times {\Omega} \right] \cap Q; R^N) \cap C_{\rm weak}([0,T]; L^2(\Omega; R^N)) ,\
\vc{v}(0, \cdot) = \vu_0 , \ \vc{v}(T, \cdot) = \vu_T;
\eF
\item
\bFormula{pr2}
\partial_t \vc{v} + \Div \tn{F} = 0,
 \ \Div \vc{v} = 0 \ \mbox{in}\ \D'((0,T) \times \Omega; R^N)
\eF
for some $\tn{F} \in C_b( \left[ (0,\tau) \cup (\tau, T) \times \Omega \right] \cap Q; R^{3 \times 3}_{{\rm sym},0} )$;
\item
\bFormula{pr3}
\frac{N}{2}\lambda_{\rm max} \left[ \frac{ (\vc{v} + h[\vc{v}] ) \otimes (\vc{v} + h[\vc{v}] ) }{r[\vc{v}]} - \tn{F} + \tn{H}[\vc{v}] \right]
 < e[\vc{v}] \vphantom{\frac{1}{2}} \ \mbox{in}\ \left[ (0,\tau)   \times \Omega \right] \cap Q,
\eF
\item
\bFormula{pr4}
\sup_{(t,x) \in Q, t > \tau + s}
\frac{N}{2}\lambda_{\rm max} \left[ \frac{ (\vc{v} + h[\vc{v}] ) \otimes (\vc{v} + h[\vc{v}] ) }{r[\vc{v}]} - \tn{F} + \tn{H}[\vc{v}] \right]
 - e[\vc{v}] < 0 \ \mbox{for any}\ 0 < s < T- \tau,
\eF
\item
\bFormula{A1}
\frac{1}{2} \intO{ \frac{ |\vc{v} + \vc{h}[\vc{v}] |^2 }{r[\vc{v}]} ({\tau}, \cdot) } = \intO{ e[\vc{v}]({\tau}, \cdot) }.
\eF

\end{itemize}
\eT

\bRemark{A1}

Unlike the subsolutions considered in the proof of Theorem \ref{TE1}, the function $\vc{v}$ satisfies (\ref{A1}) and is therefore \emph{strongly}
continuous at the point ${\tau}$ attaining the desired energy $e[\vc{v}]({\tau}, \cdot)$ in the integral sense. There is no energy jump at the time
$t = {\tau}$! Moreover, the set of such times is dense in $(0,T)$.

\eR

\bRemark{A2}

In view of hypothesis (\ref{HYP}) we have
\bFormula{CON}
C_b(Q; R^M) \subset C_{\rm loc}(0,T; L^q(\Omega; R^M)) \ \mbox{for any}\ 1 \leq q < \infty.
\eF
This observation justifies (\ref{A1}) and
will be frequently used in the proof below.

\eR

\bProof

The function $\vc{v}$ will be constructed recursively as a limit
\bFormula{A1a}
\vc{v}_k \to \vc{v} \ \mbox{in}\ C_{\rm weak}([0,T]; L^2(\Omega; R^N)) \ \mbox{for suitable}\  \vc{v}_k \in X_0.
\eF

We start by fixing the open interval $(a_0,b_0) \subset (0,T)$ in which the time $\tau$ is to be localized. As the space $X_0$ of
subsolutions is non-empty, we take
\[
\vc{v}_0 \in X_0, \ \mbox{along with the associated flux}\ \tn{F}_0.
\]

Next, we construct a sequence of functions $\vc{v}_k$, open intervals $(a_k, b_k) \subset (0,T)$, times
$\tau_k \in (a_k, b_k)$, and a decreasing sequence of positive numbers $\delta_k \searrow  0$ such that:

\begin{itemize}

\item

\bFormula{A2-}
\partial_t \vc{v}_k + \Div \tn{F}_k = 0,\ \Div \vc{v}_k = 0 \ \mbox{in}\ \D'((0,T) \times \Omega), \ \vc{v}(0) = \vu_0, \ \vc{v}(T) = \vu_T,
\eF
for a certain field $\tn{F}_k \in C(Q; R^{N \times N}_{{\rm sym},0})$,
\bFormula{A2}
\begin{array}{c}
\vc{v}_k - \vc{v}_{k-1} \in \DC(Q;R^N), \ {\rm supp} [\vc{v}_k - \vc{v}_{k-1}] \subset [(a_k, b_k) \times \Omega] \cap Q,\\ \\
\tn{F}_k - \tn{F}_{k-1} \in \DC(Q;R^{N \times N}_{{\rm sym},0}), \ {\rm supp}[\tn{F}_k - \tn{F}_{k-1}] \subset [(a_k, b_k) \times \Omega] \cap Q,
\end{array}
\eF
where
\[
0 < a_{k-1} < a_k < b_k < b_{k-1}, \ \ep_k = b_k - a_k \to 0 \ \mbox{for}\ k \to \infty;
\]

\item
\bFormula{A3a}
\sup_{t \in [0,T]} d (\vc{v}_k (t) , \vc{v}_{k-1}(t) ) < \frac{1}{2^k} ;
\eF

\item

\bFormula{A3}
\sup_{t \in (0,T) } \left| \intO{ \frac{1}{r[\vc{v}_{j}]} (\vc{v}_k - \vc{v}_{k-1} ) \cdot \vc{v}_j } \right| < \frac{1}{2^k}
\eF
for all $j=0, \dots, k-1$;

\item
there exists $\tau_k \in (a_k,b_k)$ and a positive constant $\lambda$ \emph{independent} of $k$ such that
\bFormula{A5}
\frac{1}{2} \intO{ \frac{ | \vc{v}_k + \vc{h}[\vc{v}_k] |^2 }{ r[\vc{v}_k]} (\tau_k, \cdot) }
\geq \frac{1}{2} \intO{ \frac{ | \vc{v}_{k-1} + \vc{h}[\vc{v}_{k-1}] |^2 }{ r[\vc{v}_{k-1}]} (t, \cdot) }
+ \lambda \frac{1}{\ep^2_k} \alpha_k^2
\eF
\[
> \frac{1}{2} \intO{ \frac{ | \vc{v}_{k-1} + \vc{h}[\vc{v}_{k-1}] |^2 }{ r[\vc{v}_{k-1}]} (\tau_{k-1}, \cdot) }
+ \frac{1}{2} \lambda \frac{1}{\ep^2_k} \alpha_k^2 \ \mbox{for all}\ t \in (a_k, b_k),
\]
where we have set
\[
\alpha_k = \int_{a_k}^{b_k} \intO{ \left( e[\vc{v}_{k-1}] -  \frac{1}{2} \frac{| \vc{v}_{k-1} + \vc{h}[\vc{v}_{k-1}] |^2 }{
r[\vc{v}_{k-1}]}   \right) }  \dt > 0;
\]

\item
\bFormula{A5a}
\frac{N}{2}\lambda_{\rm max} \left[ \frac{ (\vc{v}_k + h[\vc{v}_{k}] ) \otimes (\vc{v}_{k} + h[\vc{v}_{k}] ) }{r[\vc{v}_{k}]} - \tn{F}_{k} + \tn{H}[\vc{v}_{k}] \right]
 < e[\vc{v}_k]  \ \mbox{in}\ [(0, a_k]  \times \Omega] \cap Q,
\eF
\bFormula{A5b}
\frac{N}{2}\lambda_{\rm max} \left[ \frac{ (\vc{v}_k + h[\vc{v}_{k}] ) \otimes (\vc{v}_k + h[\vc{v}_{k}] ) }{r[\vc{v}_{k}]} - \tn{F}_k + \tn{H}[\vc{v}_{k}] \right]
 < e[\vc{v}_{k}] - \delta_{k} \left(1 + \frac{1}{2^k} \right)
 \eF
\[
\mbox{in}\ [(a_k, b_k)  \times \Omega] \cap Q,
\]

\bFormula{A5c}
\frac{N}{2}\lambda_{\rm max} \left[ \frac{ (\vc{v}_k + h[\vc{v}_{k}] ) \otimes (\vc{v}_k + h[\vc{v}_{k}] ) }{r[\vc{v}_{k}]} - \tn{F}_k + \tn{H}[\vc{v}_{k}] \right]
 < e[\vc{v}_{k}] - \delta_{j-1} \left(1 + \frac{1}{2^k} \right)
 \eF
\[
\mbox{in}\ [[b_{j}, b_{j-1})  \times \Omega] \cap Q, \ j=0,1, \dots, k, \ b_{-1} \equiv T.
\]

\end{itemize}

\medskip

{\bf Step 1}

It follows from the properties of $X_0$ that $\vc{v}_0$ satisfies (\ref{A2-}), along with the bounds (\ref{A5a} -- \ref{A5c}) for a certain
\[
\delta_0 = \delta_{-1} > 0.
\]

\medskip

{\bf Step 2}

Suppose we have already constructed the functions $\vc{v}_j$, along with intervals $(a_j, b_j)$, the times $\tau_j$,
and the constants $\delta_j$, for $j=0,1,\dots k - 1$ enjoying the properties (\ref{A2-} -- \ref{A5c}). Our goal is to find $\vc{v}_k$, $(a_k, b_k)$,
$\tau_k$, and $\delta_k$.

First, we fix the interval $(a_k,b_k)$. To this end, compute
\[
\alpha_k =\int_{a_k}^{b_k} \intO{ \left( e[\vc{v}_{k-1}] -  \frac{1}{2} \frac{| \vc{v}_{k-1} + \vc{h}[\vc{v}_{k-1}] |^2 }{
r[\vc{v}_{k-1}]}   \right) }  \dt.
\]
As a consequence of (\ref{CON}), the integrand is a continuous function of time continuous function of time; whence
\[
\frac{\alpha_k}{\ep_k} = \frac{1}{\ep_k} \int_{a_k}^{b_k} \intO{ \left( e[\vc{v}_{k-1}] -  \frac{1}{2} \frac{| \vc{v}_{k-1} + \vc{h}[\vc{v}_{k-1}] |^2 }{
r[\vc{v}_{k-1}]}   \right) }  \dt
\]
\[
\to
\]
\[
\intO{ \left( e[\vc{v}_{k-1}] -  \frac{1}{2} \frac{| \vc{v}_{k-1} + \vc{h}[\vc{v}_{k-1}] |^2 }{
r[\vc{v}_{k-1}]}   \right)(\tau_{k-1}) } \ \mbox{as}\ \ep_k = b_k - a_k \to 0.
\]
Consequently, keeping in mind that $\alpha_k > 0$ and repeating the same continuity argument, we may choose $a_{k-1} < a_k < b_k < b_{k-1}$ and
$\ep_k$ so small that
\bFormula{A6}
\frac{1}{\ep_k} \int_{a_k}^{b_k} \intO{ \frac{1}{2} \frac{ | \vc{v}_{k-1} + \vc{h}[\vc{v}_{k-1}] |^2}{r[\vc{v}_{k-1}]} } \dt
+ \Lambda(\Ov{e}) \frac{\alpha_k^2}{\ep^2_k}
\eF
\[
\geq \intO{ \frac{1}{2} \frac{ | \vc{v}_{k-1} + \vc{h}[\vc{v}_{k-1}] |^2}{r[\vc{v}_{k-1}]}(t, \cdot) }
+ \Lambda(\Ov{e}) \frac{\alpha_k^2}{2 \ep^2_k}
\]
\[
\geq
\intO{ \frac{1}{2} \frac{ | \vc{v}_{k-1} + \vc{h}[\vc{v}_{k-1}] |^2}{r[\vc{v}_{k-1}]}(\tau_{k-1}, \cdot) }
+ \Lambda(\Ov{e}) \frac{\alpha_k^2}{4 \ep^2_k} \ \mbox{for all}\ t \in (a_k,b_k),
\]
where $\Lambda(\Ov{e})$ is the universal constant introduced in Lemma \ref{LO1}.

At this stage, we apply Lemma \ref{LO1} for
\[
U = [(a_k, b_k) \times \Omega] \cap Q, \ \tilde{\vc{h}} = \vc{v}_{k-1} + \vc{h}[\vc{v}_{k-1}],\ \tilde {r} = r[\vc{v}_{k-1}],\ \tilde{\tn{H}} = \tn{F}_{k-1} -\tn{H}[\vc{v}_{k-1}],
\]
and
\[
\tilde e = e[\vc{v}_{k-1}] - \delta_k \left(1 + \frac{1}{2^{k-1}} \right),
\]
where $\delta_k > 0$ is chosen small enough so that (\ref{O12-}) may hold.

Now we claim that it is possible to take
\bFormula{DEF}
\vc{v}_k = \vc{v}_{k-1} + \vc{w}_n, \ \tn{F}_k = \tn{F}_{k-1} + \tn{G}_n, \ n \ \mbox{large enough},
\eF
where $\vc{w}_n$, $\tn{G}_n$ are the quantities constructed in Lemma \ref{LO1}. Obviously, the functions $\vc{v}_k$ satisfy (\ref{A2-} -- \ref{A3}) provided $n$ is large enough. Indeed we observe that $n$ can be chosen so large for (\ref{A3}) to be satisfied. To see this we realize that, by virtue of (\ref{CON}), the image
\[
\cup_{t \in [a_k,b_k]} \frac{ \vc{v}_j }{r[\vc{v}_j]} (t, \cdot) \ \mbox{is compact in}\ L^2(\Omega; R^N), \ j=0,\dots, k-1.
\]

Next, we use continuity of the operators $\vc{h}$, $r$ specified in (\ref{CON}) to compute
\[
\int_{a_k}^{b_k} \intO{ \frac{1}{2} \frac{|\vc{v}_k + \vc{h}[\vc{v}_{k}]|^2}{r [\vc{v}_{k}] } } \ \dt
= \int_{a_k}^{b_k} \intO{ \frac{1}{2} \frac{|\vc{v}_{k-1} + \vc{h}[\vc{v}_{k-1} + \vc{w}_n]|^2}{r [\vc{v}_{k-1} + \vc{w}_n] } } \ \dt
\]
\[
+ \int_{a_k}^{b_k} \intO{ \frac{1}{2} \frac{|\vc{w}_n|^2}{r [\vc{v}_{k-1} + \vc{w}_n] } } \ \dt + 2 \int_{a_k}^{b_k} \intO{
\frac{ \vc{w}_n \cdot (\vc{v}_{k-1} + \vc{h}[\vc{v}_{k-1} + \vc{w}_n])}{r[\vc{v}_{k-1} + \vc{w}_n]} } \ \dt
\]
\[
\geq \int_{a_k}^{b_k} \intO{ \frac{1}{2} \frac{|\vc{v}_{k-1} + \vc{h}[\vc{v}_{k-1}]|^2}{r [\vc{v}_{k-1}] } } \ \dt +
\frac{1}{2} \int_{a_k}^{b_k} \intO{ \frac{1}{2} \frac{|\vc{w}_n|^2}{r [\vc{v}_{k-1}] } } \ \dt + e_n
\]
provided $n$ is large enough, where $e_n \to 0$ as $n \to \infty$ for any fixed $k$. Consequently, (\ref{O1}) gives rise to
\bFormula{A7}
\int_{a_k}^{b_k} \intO{ \frac{1}{2} \frac{|\vc{v}_k + \vc{h}[\vc{v}_{k}]|^2}{r [\vc{v}_{k}] } } \ \dt + e_n \geq \int_{a_k}^{b_k} \intO{ \frac{1}{2} \frac{|\vc{v}_{k-1} + \vc{h}[\vc{v}_{k-1}]|^2}{r [\vc{v}_{k-1}] } } \ \dt
\eF
\[
+
\frac{\Lambda(\Ov{e})}{4} \int_{a_k}^{b_k} \intO{ \left( e[\vc{v}_{k-1}] - \frac{1}{2}
\frac{ |\vc{v}_{k-1} + \vc{h}[\vc{v}_{k-1}] |^2}{r[\vc{v}_{k-1}]} \right)^2 } \ \dt
\]
\[
\geq
\int_{a_k}^{b_k} \intO{ \frac{1}{2} \frac{|\vc{v}_{k-1} + \vc{h}[\vc{v}_{k-1}]|^2}{r [\vc{v}_{k-1}] } } \ \dt +
\frac{\Lambda(\Ov{e})} {4} \frac{1}{|\Omega|} \frac{\alpha_k^2}{\ep_k},
\]
where the last line follows from Jensen's inequality. Thus,
using (\ref{A6}), (\ref{A7}), we may find $n$ large enough and $\tau_k \in (a_k, b_k)$ such that (\ref{A5}) holds with some
$\lambda$ that can be determined in terms of $\Lambda(\Ov{e})$ and $|\Omega|$.

\medskip

Finally, our goal is to check that $\vc{v}_k$, $\tn{F}_k$ satisfy (\ref{A5a} -- \ref{A5c}). First we claim that (\ref{A5a}) is a direct consequence of the causality property (\ref{E4b}). Next, Lemma (\ref{LO1}), specifically (\ref{O1-}), yields
\[
\frac{N}{2} \lambda_{\rm max} \left[ \frac{ (\vc{v}_{k} + \vc{h}[ \vc{v}_{k-1} ] ) \otimes (\vc{v}_{k} + \vc{h}[ \vc{v}_{k-1} ] )}{{r}[\vc{v}_{k-1}]} - \tn{F}_k + \tn{H}[\vc{v}_{k-1}] \right] < e[\vc{v}_{k-1}]  - \delta_k \left(1 + \frac{1}{2^{k-1}} \right)
\]
in $[(a_k,b_k) \times \Omega] \cap Q$; whence (\ref{A5b}) follows from uniform continuity of $\vc{h}$, $r$, $\tn{H}$ and $e$ provided $n$ is large enough.
To see (\ref{A5c}), we have to realize that $\vc{v}_k = \vc{v}_{k-1}$ and $\tn{F}_k = \tn{F}_{k-1}$ in $[(b_k, T) \times \Omega] \cap Q$, and, similarly to the above, relation (\ref{A5c}) follows from continuity of  $\vc{h}$, $r$, $\tn{H}$ and $e$ as son as $n$ is chosen large enough.

\medskip

{\bf Step 3}

Our ultimate goal is to observe that $\vc{v}$, determined by the limit (\ref{A1a}), enjoys the desired properties claimed in Theorem \ref{TA1}. We set
$\tau = \lim_{k \to \infty} \tau_k$. Since the functions $\vc{v}_k$, $\tn{F}_k$ coincide with $\vc{v}_{k-1}$, $\tn{F}_{k-1}$ on the time intervals
$(0, a_k)$, $(b_k,T)$, the properties (\ref{pr1} -- \ref{pr4}) follow by taking the limit in (\ref{A2-}), (\ref{A5a} -- \ref{A5c}) for $k \to \infty$.

To see (\ref{A1}), we first observe that, by virtue of (\ref{A5}),
\[
\frac{1}{2} \intO{ \frac{ | \vc{v}_k + \vc{h}[\vc{v}_k ] |^2 }{r [\vc{v}_k]} (\tau_k) } \nearrow Y \ \mbox{as}\ k \to \infty,
\]
and the convergence is uniform on the time intervals $(a_k, b_k)$; whence
\[
\frac{\alpha_k}{\ep_k} = \frac{1}{\ep_k} \int_{a_k}^{b_k} \intO{ \left( e[\vc{v}_{k-1}] -  \frac{1}{2} \frac{| \vc{v}_{k-1} + \vc{h}[\vc{v}_{k-1}] |^2 }{
r[\vc{v}_{k-1}]}   \right) }  \dt \to 0,
\]
which in turn implies
\bFormula{fA1}
\frac{1}{2} \intO{ \frac{ |\vc{v}_k + \vc{h}[\vc{v}_k ]|^2 }{r [\vc{v}_k]} (t) } \to \intO{ e[\vc{v}](\tau) } \ \mbox{as}\ k \to \infty
\ \mbox{uniformly for}\ t \in (a_k, b_k).
\eF

We show that (\ref{fA1}) yields
\[
\vc{w}_k (\tau, \cdot) \to \vc{w}(\tau, \cdot)
\]
which completes the proof of Theorem \ref{TA1}. Indeed we may write
\[
\intO{ \frac{| \vc{v}_m - \vc{v}_n |^2 }{r[\vc{v}_n]} } = \intO{ \frac{| \vc{v}_m |^2 }{r[\vc{v}_n]} } - \intO{ \frac{| \vc{v}_n |^2 }{r[\vc{v}_n]} }
- 2 \intO{ \frac{( \vc{v}_m - \vc{v}_n ) \cdot \vc{v}_n }{r[\vc{v}_n]} }, \ m > n,
\]
where the difference of the first two integrals vanishes as $n \to \infty$ uniformly for $t \in (0,T)$; whereas
\[
\left| \intO{ \frac{( \vc{v}_m - \vc{v}_n ) \cdot \vc{v}_n }{r[\vc{v}_n]} } \right| = \left| \sum_{k=n}^{m - 1} \intO{ \frac{( \vc{v}_{k+1} - \vc{v}_k ) \cdot \vc{v}_n }{r[\vc{v}_n]} } \right| \leq \frac{1}{2^{n-1}} \ \mbox{uniformly in}\ (0,T)
\]
in view of (\ref{A3}).

\rightline{\qed}

Now, we can define a set of subsolutions on the time interval $(\tau, T)$, with
\[
\vu_0 = \vc{v}(\tau), \ Q = Q_\tau = Q \cap \left[(\tau, T) \times \Omega\right],
\]
and the operators $\vc{h}_\tau$, $r_\tau$, $e_\tau$, $\tn{H}_\tau$ defined as
\[
\vc{h}_\tau [ \vc{w} ] = \vc{h} [ \tilde \vc{w} ]|_{(\tau, T)}, \ \mbox{where} \ \tilde \vc{w} = \left\{
\begin{array}{l} \vc{v} \ \mbox{in} \ [0, \tau] \\ \\ \vc{w} \in [\tau, T] \end{array} \right. ,
\]
where $\vc{v}$ is the function constructed in Theorem \ref{TA1}. In accordance with (\ref{A1}), we have $\vc{v}|_{[\tau, T]}$ is
a subsolution, and
\[
\frac{1}{2} \intO{ \frac{|\vc{u}_0 + \vc{h}_\tau [\vc{u}_0] |^2 }{r_\tau [\vc{u_0}]}  } = \intO{ e_\tau[\vc{u}_0] }.
\]

Finally we note that in this case the weak solutions $\vc{u}$ constructed via Theorem \ref{TE1} will satisfy
\[
\frac{1}{2} \frac{ | \vu + \vc{h} [\vu] |^2 }{r[\vu]}  (t, x) = e[\vu] (t,x) \ \mbox{for a.a.}\ (t,x) \in (0,T) \times \Omega
\ \mbox{and \emph{including} the initial time}\ t = 0.
\]

\subsection{Example, dissipative solutions to the Euler-Fourier system}
\label{KON}

Revisiting the Euler-Fourier system introduced in Section \ref{EF}, we say that $\vr$, $\vt$, $\vu$ is a \emph{dissipative solution} of (\ref{X1} -- \ref{X3}),
if, in addition, the energy balance
\bFormula{KEX1}
E(t) \equiv \intO{ \left( \frac{1}{2} \vr |\vu|^2 + \frac{3}{2} \vr \vt \right)(t, \cdot) } =
\intO{ \left( \frac{1}{2} \vr_0 |\vu_0|^2 + \frac{3}{2} \vr_0 \vt_0 \right)}
\eF
holds for a.a. $t \in (0,T)$.

As a possible application of Theorem \ref{TA1}, one can show the following result, see \cite[Theorem 4.2]{ChiFeiKre}:

\bTheorem{KEX1}
Under the hypotheses of Theorem \ref{TX1}, let $T > 0$ and the data
\[
\vr_0, \ \vt_0 \in C^2(\Omega), \ \vr_0, \vt_0 > 0
\]
be given.

Then there exists $\vu_0 \in L^\infty(\Omega; R^N)$ such that the Euler-Fourier system (\ref{X1} -- \ref{X3}), with the initial conditions (\ref{X7}),
admits infinitely many \emph{dissipative} solutions in $(0,T) \times \Omega$.

\eT

\def\cprime{$'$} \def\ocirc#1{\ifmmode\setbox0=\hbox{$#1$}\dimen0=\ht0
  \advance\dimen0 by1pt\rlap{\hbox to\wd0{\hss\raise\dimen0
  \hbox{\hskip.2em$\scriptscriptstyle\circ$}\hss}}#1\else {\accent"17 #1}\fi}


\end{document}